\documentclass[oneside,english]{amsart}
\usepackage[LGR,T1]{fontenc}
\usepackage[latin9]{inputenc}
\usepackage{amsthm}
\usepackage{amssymb}
\usepackage{stackrel}

\makeatletter

\DeclareRobustCommand{\greektext}{%
  \fontencoding{LGR}\selectfont\def\encodingdefault{LGR}}
\DeclareRobustCommand{\textgreek}[1]{\leavevmode{\greektext #1}}

\numberwithin{equation}{section}
\numberwithin{figure}{section}

\makeatother

\usepackage{babel}
\begin{document}
\title{Estimations of Fourier Coefficients For Controlled Distortion of a
Periodic Function.}
\author{Vladimir Sluchak}

\maketitle

There is a class of physical filtration processes where the input
is adequately modeled by a continuous periodic function \textit{f
(x)} of bounded variation over its period, and the output depends
only on certain harmonics of the Fourier expansion of \textit{f (x)}
in the orthogonal basis of trigonometric functions. One example of
such a process is the discrete spectrum sound generation by a revolving
body in a steady fluid flow. This type of sound can be controlled
through the amplitudes of certain harmonics of the circular distribution
of the inflow velocity \cite{bavin1983marine}. Accordingly, to formulate
attainable goals of such a controlled distortion of \textit{f (x)},
some inequalities relating the integral characteristics of \textit{f
(x)} with the coefficients of the corresponding Fourier expansion
would be useful.

\textit{
\[
f(x)=\frac{a_{0}}{2}+\mathop{\stackrel[i=1]{\infty}{\sum}(a_{j}\cos}jx+b_{j}\sin jx);a_{j}=\frac{1}{\pi}\stackrel[0]{2\pi}{\int}f(x)\cos\begin{aligned}jx & dx\end{aligned}
;b_{j}=\frac{1}{\pi}\stackrel[0]{2\pi}{\int}f(x)\sin\begin{aligned}jx & dx\end{aligned}
\]
}
\begin{flushleft}
Riemann-Lebesgue Theorem \cite{rudin1987real} gives the order \textit{O
(1/j)} for the amplitudes of the harmonics $a_{j}$ and $b_{j}$ as
$j$\ensuremath{\in}$\mathbb{N}$, $j\rightarrow\infty$. 
\par\end{flushleft}

Taking into consideration this theorem, we will seek to prove the
inequalities of the following structure:

\begin{equation}
|a_{j}|\leqslant\frac{F(f(x))}{j}
\end{equation}

where \textit{F (f (x))} is a \emph{Functional} of the integral characteristics
of \textit{f (x)} on the interval of expansion. 

Breaking down the interval of expansion {[}0,2\textgreek{\textpi}{]}
into segments {[}$\frac{\pi(k-1)}{2j}$,$\frac{\pi k}{2j}${]}, k=1,2\dots 4j,
where functions \textit{sin (jx)} and \textit{cos (jx)} do not change
sign we get (for cosine):
\begin{flushleft}
\begin{equation}
a_{j}=\frac{1}{\pi}\sum_{k=1}^{4j}\int\limits_{\frac{\pi(k-1)}{2j}}^{\frac{\pi{k}}{2j}}f(x)\cos\ jx\ dx
\end{equation}
\par\end{flushleft}

And using generalized mean value theorem for definite integrals:

\[
a_{j}=\frac{1}{\pi}\stackrel[k=1]{4j}{\sum}f(x_{k})\stackrel[\frac{\pi(k-1)}{2j}]{\frac{\pi k}{2j}}{\int}\cos\begin{aligned}jx & dx=\frac{1}{\pi j}\stackrel[k=1]{4j}{\sum}f(x_{k})\left[\sin\frac{\pi k}{2}-\sin\frac{\pi\left(k-1\right)}{2}\right]\end{aligned}
\]

\begin{center}
\begin{equation}
a_{j}=\frac{1}{\pi j}\stackrel[k=1]{4j}{\sum}\left(-1\right)^{k\setminus2}f(x_{k})
\end{equation}
 where $\frac{\pi(k-1)}{2j}\leqslant x_{k}$$\leqslant\frac{\pi k}{2j}$ 
\par\end{center}

Formula (3) can be presented as 

\begin{equation}
a_{j}=\frac{1}{\pi j}\stackrel[k=1]{2j}{\sum}(-1)^{k-1}\left[f(x_{2k-1})-f(x_{2k})\right]
\end{equation}

\begin{center}
where $\left[x_{2k-1};x_{2k}\right]$$\in$$\left[\frac{\pi(k-1)}{j};\frac{\pi k}{j}\right]$
\par\end{center}

\begin{flushleft}
As modulus of sum is less then or equal to the sum of modulus the
following inequality derives from (4):
\par\end{flushleft}

\begin{equation}
|a_{j}|\leq\frac{1}{\pi j}\sum_{k=1}^{2j}|f(x_{2k-1})-f(x_{2k})|
\end{equation}

\begin{center}
where $\left[x_{2k-1};x_{2k}\right]\in\left[\frac{\pi(k-1)}{j};\frac{\pi k}{j}\right]$
\par\end{center}

As the segments $\left[x_{2k-1};x_{2k}\right]$ are not overlapping,
the sum $\sum_{k=1}^{2j}|f(x_{2k-1})-f(x_{2k})|$ is less than $V_{0}^{2\pi}(f)$
- the full variation of $f\left(x\right)$ on the interval $[0,2\pi]$.
Hence the following inequality:

\begin{equation}
|a_{j}|\leq\frac{V_{0}^{2\pi}(f)}{\pi j}
\end{equation}

where $\mathop{V_{0}^{2\pi}(f)}$=$\stackrel[0]{2\pi}{\int}|f^{\prime}(x)|dx$
and $\pi$ is the square of the\textit{ norm }for the orthogonal basis
\textit{sin (jx), cos (jx)}

Note that full variation of $f(x)$ is equal to the sum of absolute
differences between closest extrema, so (6) can be presented as:

\begin{equation}
|a_{j}|\leq\frac{1}{\pi j}\sum_{k=1}^{N}\Delta_{k}
\end{equation}

where $\Delta_{k}$is a difference between closest extrema, N - total
number of extrema of $f(x)$ on the interval of expansion $\left[0,2\pi\right].$

One more weak, but not trivial restriction follows from (5):

\begin{equation}
|a_{j}|\leq\frac{2}{\pi}\Delta
\end{equation}

where $\Delta=maxf(x)-minf(x),x\in\left[0,2\pi\right]$.

It is easy to show that the inequalities (6), (7) and (8) derived
for the cosine amplitudes are valid for the sines.

These inequalities allow to build fuzzy targets for the controlled
distortion of $f(x)$ when the control process is aimed at reducing
the amplitudes of certain harmonics of the corresponding Fourier expansion.

Given $f(x)$ and its Fourier expansion as described, let the aim
of the controlled distortion be - to reduce the amplitude $|a_{j}|$.
Further we denote the target function - $\tilde{f}(x)$.

The ideal result of such a controlled distortion would be $a_{j}=0$,
which implies a class of the continuous functions with bounded variation
$f^{a_{j}=0}(x)$ and the difference $f_{\triangle}(x)=\tilde{f}(x)-f^{a_{j}=0}(x)$. 

Note, that $f^{a_{j}=0}(x)$, depending on the control process, can
vary from \textit{``trivial''} \textit{- const} on the interval
$x\in\left[0,2\pi\right]$ to $f(x)-a_{j}cos(jx)$ - ''\textit{minimal}''
in terms of squared deviation from the initial $f(x)$.

Let the aim of the controlled distortion of $f(x)$ be - to reduce
$|a_{j}|$ by no less than q times, then from (7) it will be enough
for the full variation of $f_{\triangle}(x)$ to satisfy the following
condition

\begin{equation}
V_{0}^{2\pi}(f_{\triangle})\leq|a_{j}^{0}|\frac{\pi j}{q}
\end{equation}

Where $|a_{j}^{0}$| - initial value of $|a_{j}$| in the Fourier
expansion of $f(x)$.

Consider a proximity band for $f^{a_{j}=0}(x)$ with equidistant boundaries
$f^{a_{j}=0}(x)-\delta/2;f^{a_{j}=0}(x)+\delta/2;\delta=const$ and
let the distorted function $\tilde{f}(x)$ fit completely inside the
band, then the width of the band satisfying (9) and so \textit{- sufficient}
for the required reduction can be derived as follows:
\begin{equation}
V_{0}^{2\pi}(f_{\triangle})\leq N\delta\leq|a_{j}^{0}|\frac{\pi j}{q}
\end{equation}

\begin{equation}
\delta=|a_{j}^{0}|\frac{\pi j}{qN}
\end{equation}

where $N$ - number of extrema of $f_{\triangle}(x)$, $q$- the target
factor of reduction for $|a_{j}|$ . Note that in case of $N=0$ the
target function $\tilde{f}(x)$ is equidistant to $f^{a_{j}=0}(x)$
on $x\in\left[0,2\pi\right]$ and so they

have identical Fourier coefficients $a_{j}$for $j$\ensuremath{\in}$\mathbb{N}$.

From (8) we can get another formula for $\delta$:
\begin{equation}
\delta=|a_{j}^{0}|\frac{\pi}{2q}
\end{equation}

So the bandwidth $\delta$ for the distortion of $f(x)$ \textit{sufficient}
for the reduction of a harmonic formulated as $\frac{|a_{j}^{0}|}{q}$
can be found from (11) or (12), where the band is 

built as \{$f^{a_{j}=0}(x)-\delta/2;f^{a_{j}=0}(x)+\delta/2$\}. 

Note that the larger value of $\delta$ requires less deviation from
the initial $f(x)$ and so (depending on the method of control) can
be more energy efficient. 

This approach to the controlled distortion of a circular distribution
of the inflow velocity was used for the reduction of sound generated
by rotating bodies in the fluid flow \cite{sluchak1986method}.

Inequality (6) can be generalized for the expansions of $f(x)$ into
the series by other orthogonal bases, at least for those that allow
for the proof path(2)$\rightarrow$(6). 

Consider the expansion of a continuous function$f(x)$ on $x\in[-1,1]$
into a series by Chebyshev polynomials of the 1st kind $T_{j}(x)=cos(j*arccos(x))$:

\begin{equation}
f(x)=\frac{a_{0}}{2}+\mathop{\stackrel[j=1]{\infty}{\sum}a_{j}T_{j}(x)};a_{j}=\frac{2}{\pi}\stackrel[-1]{1}{\int}\frac{f(x)cos(j*arccos(x))}{\sqrt{1-x^{2}}}dx
\end{equation}

After substitution $x=cos\theta$ the coefficients of the expansion
can be presented as:
\begin{equation}
a_{j}=\frac{2}{\pi}\stackrel[0]{\pi}{\int}f(cos\theta)cos(j\theta)d\theta
\end{equation}

Then making same transformation as in (2), (3), we receive:

\begin{equation}
a_{j}=\frac{2}{\pi j}\stackrel[k=1]{2j}{\sum}\left(-1\right)^{k\setminus2}f(cos(\theta_{k}))
\end{equation}

\begin{center}
where $\frac{\pi(k-1)}{2j}\leqslant\theta_{k}\leqslant\frac{\pi k}{2j}$
\par\end{center}

And similarly to (4) and (5):
\begin{equation}
|a_{j}|\leq\frac{2}{\pi j}\sum_{k=1}^{j}|f(cos(\theta_{2k-1}))-f(cos(\theta{}_{2k}))|
\end{equation}

\begin{center}
where $\left[\theta_{2k-1};\theta_{2k}\right]\in\left[\frac{\pi(k-1)}{j};\frac{\pi k}{j}\right]$
\par\end{center}

As the segments $\left[\theta_{2k-1};\theta_{2k}\right]$ are not
overlapping, the sum $\sum_{k=1}^{j}|f(cos(\theta_{2k-1}))-f(cos(\theta{}_{2k}))|$
is less than $V_{0}^{\pi}(f(cos\theta))$ - the full variation of
$f\left(cos\theta\right)$ on the interval $[0,\pi]$.

Full variation $V_{0}^{\pi}(f(cos\theta))$ is equal to $V_{-1}^{1}(f(x))$:

\begin{equation}
V_{0}^{\pi}(f(cos\theta))=\stackrel[0]{\pi}{\int}|f^{'}(cos\theta)|d\theta=\stackrel[0]{\pi}{\int}|\frac{df(cos\theta)}{d(cos\theta)}|d(cos\theta)=\stackrel[-1]{1}{\int}|df(x)|=V_{-1}^{1}(f(x))
\end{equation}

This way, for the coefficients of expansion of a continuous periodic
function $f(x)$ into Fourier series by Chebyshev polynomials of the
1st kind, we get generally the same inequality as for the case of
the expansion into a trigonometric basis (6). 
\begin{equation}
|a_{j}|\leq\frac{2V_{-1}^{1}(f(x))}{\pi j}
\end{equation}

And in more generic form:

\begin{equation}
|a_{j}|\leq\frac{V^{ort}(f(x))}{j\parallel N\parallel^{2}}
\end{equation}

where $V^{ort}(f(x))$ is a full variation of function $f(x)$ on
the interval of orthogonality of the basis, ||N|| - norm of the basis,
j $\geqslant1$ - the number of the coefficient (defined by the number
of zeros of the member of the orthonormal basis, on the interval of
orthogonality). 

\bibliographystyle{plain}
\nocite{*}
\bibliography{books}

\end{document}